\documentclass[12pt]{amsart}
\begin{document}

\newcommand{\de}{D_{E}}
\newcommand{\df}{D_{F}}
\newcommand{\dpush}{\vec{D}}
\newcommand{\pf}{\phi(D_{E} \oplus D_{I^{\perp}})}
\newcommand{\po}{\phi(D_{F})}
\newcommand{\ho}{Hom^{\times} (E,F)}
\newcommand{\uni}{\phi (D_{\pi^{\ast} F})}
\newcommand{\unf}{\phi(D_{\pi^{\ast} E} \oplus D_{\tilde{I}^{\perp}})}
\newcommand{\res}{Res_{\phi,i}}
\newcommand{\sig}{\Sigma_{i}}
\newcommand{\stief}{V_{n,m}}
\newcommand{\gras}{G_{n,m}}
\newcommand{\cut}{X-\Sigma}
\newcommand{\cuts}{X - \bigcup \Sigma_{i}}
\newcommand{\pis}{\pi_{i \ast}}                                                                                                                                                                                                                                                                                                                                                                                                                                                                                                                                        \newcommand{\su}{\sum_{i=1}^{l}}
\newcommand{\tb}{N_{i,\epsilon}}
\newcommand{\bo}{\partial N_{i,\epsilon}}
\newcommand{\pfs}{\phi(\pi^{\ast} D_{\tilde{E}})}
\newcommand{\pos}{\phi(\pi^{\ast} D_{\tilde{F}})}
\newcommand{\rest}{\tilde{Res}_{\phi,i}}
\newcommand{\loc}{L_{loc}^{1}}
\newcommand{\im}{\lim_{\epsilon \rightarrow 0}}
\newcommand{\norm}{\parallel v \parallel}
\newcommand{\on}{\quad \mbox{on} \quad}
\newcommand{\rb}{\mathbb R}
\newcommand{\cb}{\mathbb C}
\renewcommand\theequation{\arabic{section}.\arabic{equation}}

\theoremstyle{plain}
\newtheorem{theorem}{Theorem}[section]
\newtheorem{corollary}[theorem]{Corollary}
\newtheorem{lemma}[theorem]{Lemma}

\theoremstyle{definition}
\newtheorem{remark}[theorem]{Remark}
\newtheorem{defn}[theorem]{Definition}
\title{Geometric Residue Theorems for Bundle Maps}
\author{Sunil Nair}
\address{Mathematics Section\\ I.C.T.P.\\P.O. Box 586\\34100 Trieste\\Italy}
\email{sunil@ictp.trieste.it}

\maketitle
\begin{abstract}
In this paper we prove geometric residue theorems for bundle maps over a compact manifold. The theory developed associates residues to the singularity submanifolds of the map for any invariant polynomial. The theory is then applied to a variety of settings: smooth maps between equidimensional manifolds, CR-singularities, finite singularities and singularities of odd forms as spinor bundle maps.
\end{abstract}
\nopagebreak

\section{Introduction}
 The prototype of residue theorems in geometry is the classical theorem of Hopf's [H]  which relates the zeroes of a vector field on a compact manifold to its Euler characteristic. In general, residue theorems associate topological invariants to the singularities of geometric objects. The aim of this paper is
to study singularities of maps between bundles over an oriented manifold, in particular, to obtain residue theorems for such singularities.

The general setting is as follows: Let $E$ and $F$ be vector bundles, real or
complex, over a compact, oriented manifold $X$, and let $\alpha:E \mapsto F$ be a bundle
map between $E$ and $F$ that drops rank outside a closed submanifold $\Sigma$. The
theory developed in this paper uses the notion of a pushforward connection, as developed by Harvey-Lawson [HL1],  to compare the characteristic classes of $E$ and $F$ and relate them to $\Sigma$. When rank(E) = rank(F) and $\phi$ is an
Ad-invariant polynomial, formulae of the type
  \[\phi(\Omega_{F}) - \phi(\Omega_{E}) = Res_{\phi}[\Sigma] + dT \]
are derived, where: $[\Sigma]$ is the current associated to the submanifold $\Sigma$,  $Res_{\phi}$ is a closed current computed in the curvatures $\Omega_{E}$ and $\Omega_{F}$ of $E$ and $F$ and the twisting of the map $\alpha$, and where $T$ is a canonical  transgression form. 

The theory of generic bundle maps has been studied in great detail by Harvey-Lawson [HL1], [HL2] and by Harvey-Semmes [HS]. This paper uses several key ideas developed in the papers above, namely the notion of pushforward connections and the universal setting for residue theorems, but takes a different point of view. While the authors above study atomic maps between bundles, for which the singularities are of the expected codimension, we allow the singularities to be nongeneric, asking only that they be closed submanifolds. This is not a strong restriction because a minor modification of the standard Thom Transversality Theorem (see [HL2]) shows that the set of smooth bundle maps which vanish nondegenerately is open and dense in the $C^{1}$-topology.

In this general setting one cannot hope to obtain residue theorems for bundle maps but the key point of this paper is that up to homotopy it is always possible to do so. More precisely, up to homotopy of the bundle map $\alpha$, and the connections $D_{E}$ and $D_{F}$ of the bundles $E$ and $F$, we can write residue formulae for any Ad-invariant polynomial $\phi$.

This approach leads to many interesting formulae and applications. To name a few, we obtain a generalized Hopf index theorem for bundle maps with finite singularities, a generalized Riemann-Hurwitz formula for smooth maps between manifolds of the same dimension, residue formulae for CR-singularities and residue formulae for Clifford and Spin bundles.

The author wishes to thank Blaine Lawson for introducing him to the subject and for his invaluable help in shaping the results obtained in his thesis, from which this paper stems. Also, he wishes to thank ICTP, where a part of this work was realized, for its hospitality.

\section{Pushforward Connections}
\setcounter{equation}{0}
In this section we review the concept of a pushforward connection introduced in [HL1]. The material covered in \S 1 to \S 3 can be found in the above paper in much greater detail.

A notational convention: Throughout this paper $X$ will denote  a manifold which is oriented and $\Sigma$ a submanifold of $X$.

Suppose $E$ and $F$ are vector bundles (real or complex) over an oriented manifold          $X$. Let $\de$, respectively $\df$, be a smooth connection on $E$, respectively $F$. Let $\alpha$ be a bundle map between $E$ and $F$ that is injective outside a submanifold $\Sigma$  of $X$. Then we can transplant the connection $\de$ to define a pushforward connection on $F$ outside the singularity set $\Sigma$  as follows.
\begin{equation}
  \dpush = \alpha \de \beta  +  \df (1 - \alpha \beta)
\end{equation}
Here $\beta$ is the `inverse of $\alpha$'. This is made precise below. Suppose $E$ and $F$  are equipped with metrics, not necessarily compatible with the connections $\de$ and $\df$. On the complement of $\Sigma$  let $I=im\alpha$ denote the image subbundle of $F$. We can now choose $\beta$ to be the orthogonal projection of $F$ onto $I$ followed by the inverse of the map
 \[\alpha : E \mapsto I\]
The transplanted connection $\dpush$  is singular because the map $\beta$ is singular on $\Sigma$. A more concrete formula for $\beta$ is given by  
  $$\beta=(\alpha^{\ast} \alpha)^{-1} \alpha^{\ast}$$  
where $\alpha^{\ast}$ denotes the adjoint of $\alpha$.

For the most part, this paper will concentrate on the equirank case, i.e, when rank(E) = rank(F). Then the singular pushforward connection $\dpush$ on $F$ is given by the simple formula 
 $$\dpush = \alpha \de \alpha^{-1} \on  \cut$$

On $\cut$ the pushforward connection $\dpush$  can be written in block form with respect to the splitting $F = I \oplus I^{\perp}$.  The matrix form of $\dpush$
blocks as an upper triangular matrix with diagonal terms $\alpha \de \beta$ and
$(1 - P) \df (1 - P)$. Here $P$ denotes the orthogonal projection of $F$ onto $I$. Therefore for any Ad-invariant polynomial $\phi$, on the Lie algebra $g\ell_{n}(\rb)$, or $g\ell_{n}(\cb)$, we have 
 $$\phi (\dpush) = \phi (\alpha \de \beta \oplus (1-P) \df (1-P) ) \on  \cut$$
where to simplify notation we write
 $$\phi (D) \equiv \phi (\Omega)$$
This denotes the invariant polynomial evaluated on the appropriate curvature. This  will be the notational covention  adopted throughout the paper. 

When restricted to sections of $I$, $\alpha \de \beta = \alpha \de \alpha^{-1}$ is gauge equivalent to $\de$ and this implies that
\begin{equation}
 \phi (\dpush) = \pf \on  \cut
\end{equation}
where $D_{I^{\perp}} = (1 - P) \df (1 - P)$ is the connection induced on $I^{\perp} \subset F$ by $\df$.

\section{ Families of Pushforward Connections and Transgressions}
\setcounter{equation}{0}
To obtain a transgression formula via Chern-Weil theory we want to introduce a family of connections on $F$. There is a nice way of doing this by using the notion of an approximate one. By an approximate one we mean a function $\chi$ which satisfies the following properties: 
 $$\chi : [ 0 , \infty ] \mapsto [ 0 , 1 ]$$                                 which is $C^{\infty}$ on [0 , $\infty$] and satisfies

 $$\chi (0) = 0 , \chi (\infty) = 1$$
and                                                                            $$\chi' \geq 0$$
 
Given a bundle map $\alpha$ we can define approximations to the inverse of $\alpha$ based on $\chi$ by setting
  $$\beta_{s} =  \frac{\alpha^{\ast}}{s^{2}} \rho (\frac{\alpha \alpha^{\ast}}{s^{2}})$$

The family of bundle maps $\beta_{s}$ is smooth for $0 < s \leq + \infty$ on $X$ with $\beta_{\infty} = 0$ and $\beta_{0} = \beta$ on $\cut$. We can now define a family of smooth connections $\dpush_{s}$  over the entire manifold $X$ including $\Sigma $ for $0 < s \leq + \infty$ by
  $$\dpush_{s} = \alpha \de \beta_{s} + \df (1 - \alpha \beta_{s})$$
Note that $\dpush_{\infty} = \df$ and $\dpush_{0} = \dpush$ on $\cut$.

We then have a family of curvature forms $\Omega_{s}$  corresponding to the family of connections $\dpush_{s}$. 
Using standard Chern - Weil theory (see [BoC]) we can write 
 $$\phi (D_{\infty}) - \phi (\dpush) = dT \on \cut$$
and using (2.2) we can rewrite this as
\begin{equation}
\po - \pf = dT \on \cut
\end{equation}
Here T denotes the transgression form for this family of connections. The explicit form of $T$ is given as follows.
 $$T = \int_{0}^{\infty} \phi (\dot{D}_{t} ; \Omega_{t}) dt$$
where $\phi (\dot{D}_{t} ; \Omega_{t}) = \frac{d}{ds} \phi (\Omega_{t} + sD_{t}) \mid_{s = 0}$ is the complete polarization of $\phi$. The aim of this paper is to extend equation (3.1) across the singularity set $\Sigma$ to the entire manifold $X$,  thereby obtaining residue formulae.

\begin{remark}
We can also consider the case where $\alpha$ is a surjective map outside $\Sigma$. Here we can define a pullback connection $ \stackrel{\leftarrow}{D}$ on $E$ by 
$$  \stackrel{\leftarrow}{D}  = \beta\df\alpha + (1 - \beta\alpha)\de$$
where $\beta$ denotes orthogonal projection onto $T$ followed by the inverse of the map $\alpha: K^{\perp} \mapsto T$ and $K \equiv ker \alpha$ is the kernel subbundle of $E$. Again by considering families of connections, we obtain, for any invariant polynomial $\phi$,
\begin{equation}
 \phi(\de) - \phi(\df \oplus D_{K}) = dT \on \cut
\end{equation}
For the sake of clarity we omit mentioning the surjective case explicitly in the exposition that follows. The formulae are the same in both cases; the reader just has to replace (3.1) with (3.2) to obtain the result for the surjective case.
\end{remark}

\section{The Universal Setting}
\setcounter{equation}{0}
It is often useful to consider the above setting universally. By this we mean transplanting the given data to $\ho$, the bundle of injective maps from $E$ to $F$.

Let $\pi : \ho \mapsto X$ be the projection map onto the manifold $X$. Then we can pull back the bundles $E$ and $F$ by $\pi$ to obtain the bundles $\pi^{\ast} E$ and $\pi^{\ast} F$ over $\ho$. There is a tautological bundle map 
 $$\tilde{\alpha} : \pi^{\ast} E \mapsto \pi^{\ast} F$$
which at a point $\alpha \in Hom^{\times} (E_{x} , F_{x})$ above $x \in X$ is simply defined to be $\alpha$. This tautological map is injective everywhere. We can pull the connections on $E$ and $F$ back to $\pi^{\ast} E$ and $\pi^{\ast} F$ and apply Chern - Weil theory to this setting as we did in the previous section. We then have the following universal formula.
\begin{equation}
 \uni - \unf = d\tilde{T} \on \ho
\end{equation}

A smooth bundle map $\alpha : E \mapsto F$ which is injective outside $\Sigma \subset X$ defines a cross-section of $\ho$ on $\cut$ and we have that
$$\alpha^{\ast} (\pi^{\ast} E) = E$$                                                          and                                                                     $$\alpha^{\ast} (\pi^{\ast}F) = F$$                                                              over $\cut$
as bundles with connections. Furthermore
 $$\alpha^{\ast} (\tilde{\alpha}) = \alpha \on \cut$$
Thus every case is a pullback of the universal one, in particular, (4.1) pulls back to give (3.1) on $\cut$.

\section{Residue Formulae}
\setcounter{equation}{0}
In this section we study equation (3.1) in more detail. We show how to obtain residue formulae when both the transgression form $T$ and the characteristic form $\po - \pf$ in the equation above extend as $\loc$ forms across the singularities of the bundle map $\alpha$.

Suppose that we are in the setting outlined in the previous chapter, where $\alpha : E \mapsto F$ is a bundle map, defined and injective outside $\bigcup\sig$. Each  $\sig$ is assumed to be a submanifold of $X$, disjoint from the others, but not necessarily of the same dimension. The submanifolds $\sig$ will be referred to as the $\bf{singularities}$ of the map $\alpha$. Then  as in    \S 3 we can write down the following:
  $$\po - \pf = dT \on \cuts$$
Suppose that the transgression form $T$   and the characteristic form $\po - \pf$ in the equation above extend as an $\loc$ forms across the singularities. We then have the following residue theorem.
\begin{theorem}
Let $\alpha : E \mapsto F$ be a map which is injective outside $\bigcup \sig$, where each $\sig$ is a closed submanifold of a compact, oriented manifold $X$.  Suppose that $T$ and $\po -\pf$ extend as  $\loc$ forms on $X$ for a given invariant polynomial $\phi$. Furthermore, assume that the extension of $\po - \pf$ is d-closed on $X$. Then
\begin{equation}
  \po - \pf = \su \res [\sig] + dT \on X
\end{equation}
where                                                                      $$\res\equiv   \im \int_{\pi_{i}\mid_{\bo}} T$$
is a closed current supported on $\sig$ and $deg(\res) = 2deg(\phi) - codim(\sig)$. Here $\bo$ denotes the boundary of an $\epsilon$-tubular neighborhood $\tb$  of $\sig$ and $\pi_{i} : \partial N_{i} \rightarrow \sig$ is projection.
\end{theorem}
\begin{proof}
  Choose $\epsilon$-tubular neigborhoods $\tb$  of $\sig$. Write
$X = ( X - \bigcup \tb ) \cup \bigcup \tb$. Since $\po - \pf$ extends as an $\loc$ form on $X$, we have 
$${\im} (\po - \pf) \wedge [X - \bigcup \tb] = (\po -\pf) \wedge [X]$$
Then we can write
\begin{eqnarray*}
 ( \po - \pf ) \wedge [X] & = & \im dT \wedge [X - \bigcup \tb] \\
                          & = &  \im d(T \wedge [X - \bigcup \tb ])  \\
                          &   & + \im \su T \wedge [  \bo ]  \\                                                           \end{eqnarray*}
Here $[X]$, $[ X - \bigcup \tb ]$ and $[  \bo ]$ denote the currents associated with  $X$, $X - \bigcup \tb$ and $ \bo$ respectively.

Since  $T$ extends as an $\loc$ form on $X$, the family $T \wedge [X - \bigcup \tb]$ converges to $T$ extended by zero, as currents on $X$. We now use the following convergence theorem found in Federer [F;4.1.19].

If $ a_{\epsilon} \rightarrow a$  and $b_{\epsilon} \rightarrow b$ as $\loc$ forms, then $da_{\epsilon} + b_{\epsilon} \rightarrow da + b$ in flat norm.

Applied here, this implies that $L_{i} \equiv \displaystyle{\im} T \wedge [\bo]$ exists in flat currents on $X$.

Since $supp( L_{i} ) \subset \sig$, this current is intrinsic to $\sig$, by Federer's flat support theorem found in Federer [F;4.1.15], i.e , $\pi_{i \ast} L_{i} = L_{i}$ where $\pi_{i} : \tb \rightarrow \sig$ is the fibration of the tubular neighborhhod over $\sig$. Now
\begin{eqnarray*}
 L_{i} & = & \pis L_{i}  \\
       & = & \pis [ \im T \wedge \bo ]  \\
       & = & \im \pis [ T \wedge \bo ]  \\
       & = & \im \int_{\pi_{i}\mid_{\bo}} T \wedge [\sig] \\
       & = & \res [\sig] \\
\end{eqnarray*}
We note that if $\sig$ is nonorientable then integration over the fiber defines a current with twisted coefficients in the orientation bundle of $\tb$. Since $X$ is orientable, $\res$ and $[\sig]$ lie in the same orientation class, hence $\res[\sig]$ is well-defined as a current on $X$.

If $\po - \pf$ extends as a closed $\loc$ form and each $\sig$ is a closed submanifold then applying the exterior derivative on both sides of (5.1) gives that each $\res$ is a closed current.
\end{proof}

\begin{remark}
We now explain integration over the fibers a little more carefully. If $\sig$ is orientable, integration over the fibers gives a closed current $\res$ on $X$ supported on $\sig$, which defines an element of $H^{\ast} (X;\rb)$. If $\sig$ is nonorientable, integration over the fibers gives a closed current on $\sig$ with twisted coefficients in the orientation bundle $o(\tb)$ of $\tb$. Taking the limit as $\epsilon \rightarrow 0$ does not change the orientation bundle $o(\tb)$, so in fact the limit is well defined. Furthermore, since 
$$ \tb \oplus T\sig = TX\mid_{\sig}$$
and $X$ is orientable, we have that 
$$o(\tb) \otimes o(T\sig) = \underline{\rb}$$
where $o(\tb)$ and $o(T\sig)$ are the orientation bundles of $\tb$ and $\sig$ respectively, and $\underline{\rb}$ is the trivial bundle. Hence $\sig$ is in the same orientation class as $\res$, so the pairing $\res[\sig]$ defines a closed current on $\sig$. Here $[\sig]$ denotes the current associated with $\sig$ with twisted coefficients in $o(\sig)$.

With this in hand, it is clear that when we say currents and differential forms, it is to be understood that we mean currents and differential forms with {\bf twisted coefficients} when $\sig$ is {\bf nonorientable}.
\end{remark}

The hypotheses of the theorem recur again throughout the paper, so for convenience we make the following definition.
\begin{defn}
We say that a bundle map $\alpha$ with singularities $\sig$ which are closed submanifolds is {\bf extendable}  if for any invariant polynomial $\phi$ the smooth forms  $T$ and $\po -\pf$ in $\cuts$ extend as  $\loc$ forms on the manifold $X$.
\end{defn}

\begin{remark} 
 If rank(E) = rank(F), then $I^{\perp}$ does not appear in the equation above. Since $\de$ and $\df$ are smooth connections, $\phi(\df) - \phi(\de)$ extends as a closed smooth form on $X$. Therefore, the only hypothesis needed in this case is that $T$ extend to be $\loc$ on $X$.
\end{remark}

\section{Normalized Maps and Normalized Bundles}
\setcounter{equation}{0}
In this section we discuss when the transgression form $T$  and the characteristic form $\po - \pf$ extend as  $\loc$ forms on $X$. This involves the notions of normalized bundles and normalized maps. We define these notions below.

First we define a tubular neighborhood structure of a closed submanifold $\Sigma$ to be an $\epsilon$-tubular neighborhood $N_{\epsilon}$ of $\Sigma$ with a given smooth identification with the normal disk bundle to $\Sigma$. Let $\pi : N \rightarrow \Sigma$ denote the bundle projection and $\rho:N-\Sigma \rightarrow \partial N$ denote the radial projection onto the boundary induced from the vector bundle structure of $\pi : N \rightarrow \Sigma$.
  
\begin{defn}
We say that a bundle $E$  with connection $\de$  over a manifold $X$ is {\bf normalized} at a submanifold $\Sigma$, if for some tubular neighborhood structure $N$  of $\Sigma$, the pair $(E,\de)$  can be written as a pullback from $\Sigma$ , i.e.,  
 $$ (E , \de) \mid_{N} \quad \equiv \quad \pi^{\ast} (E , \de) \mid_{\Sigma}$$
where $\pi : N \rightarrow \Sigma$ is  projection.
\end{defn}

Since $\Sigma$ is the retract of $N$, any bundle over $N$ is equivalent to a pullback of a bundle on $\Sigma$. However, the pullback connection is only homotopic to the original one. The condition that a bundle be normalized at $\Sigma$ ensures that the pullback connection is equal to the original one. If $\Sigma$ is a point on $X$ then the normalization condition implies that $E$ is flat in a neighborhood of the point. In general, normalization can be viewed as `flatness' in radial directions.

We now go on to define normalized maps.
\begin{defn}
A bundle map $\alpha : E \rightarrow F$ defined over $\cut$ is {\bf normalized} at $\Sigma$, a submanifold of $X$, if there exists a tubular neighborhood structure $N$ for which $(E,\de)$ and $(F,\df)$ are normalized  and such that $\alpha$ is radially constant, i.e.,
 $$ (\rho^{\ast}\alpha : E \mid_{\partial N} \rightarrow F \mid_{\partial N}) = (\alpha : E \rightarrow F) \mbox{ in } N - \Sigma$$
where $\rho : N - \Sigma \rightarrow \partial N$ is radial projection onto the boundary.
\end{defn}

With these definitions in hand we prove the main extension theorem.

\begin{theorem}
Let $\alpha : E \rightarrow F$ be an injective map outside $\bigcup \sig$. Suppose that $\alpha$ is normalized at each $\sig$. Then for any invariant polynomial $\phi$, the forms $T$ and $\po - \pf$ extend as  $\loc$ forms over the manifold $X$, i.e., $\alpha$ is extendable.
\end{theorem}
\begin{proof}
 Since $\alpha$ is normalized at each $\sig$, the family of pushforward connections $\vec{D}_{s}$  is also a pullback from $\partial \tb$, i.e.,
 $$(F, \vec{D}_{s}) \mid_{\tb - \sig} = \rho^{\ast} (F , \vec{D}_{s}) \mid_{\bo}$$
Here $\tb$ denotes the $\epsilon$-disk bundle of $\sig$ and $\bo$ its boundary.

This immediately implies that the transgression form is also a pullback from $\bo$, i.e.,
 $$ T \mid_{\tb - \sig} = \rho^{\ast} (T \mid_{\bo})$$

Now we want to show that this pullback property implies that $T$ extends as an $\loc$ form. Since this is a local property we need only construct a local argument.  Without loss of generality we can assume that each $\sig$ is of dimension 0.

Let $f : R^{n} - \{0\} \rightarrow S^{n - 1}$ be radial projection. In coordinates $f(x) = \frac{x}{\parallel x \parallel}$ where $x = (x_{i}, . . . , x_{n})$. Then
\begin{eqnarray*}
f^{\ast} dx^{i} & = & d(\frac{x_{i}}{\parallel x \parallel})\\
                & = & \frac{dx_{i}}{\parallel x \parallel} - \sum \frac{x_{i} x_{j}}{\parallel x \parallel^{2}} dx_{j}\\
\end{eqnarray*}

This is a homogeneous form of degree 0 and the coefficients of $f^{\ast} dx^{i}$ are bounded by $\frac{n}{\parallel x \parallel}$. Therefore the coefficients of $f^{\ast} dx^{I}$, for a multiindex $I$ are bounded by $(\frac{n}{\parallel x \parallel})^{\mid I \mid}$. Hence it is $\loc$ if $\mid I \mid \leq n -1$. For any $\varphi \in E^{p} (S^{n - 1})$ where
  $$ \varphi = \sum a_{I} dx^{I} \mid_{S^{n - 1}},\quad \mid I \mid \leq n -1$$
the coefficients of $f^{\ast} \varphi$ are bounded by $\sum sup  \mid a_{I} \mid    \frac{c}{\parallel x \parallel^{n -1}}$, where $c$ is a constant. Here the sup is taken over the sphere. Hence $f^{\ast} \varphi$ is $\loc$.

Applied to the transgression form  this implies that $T$ extends as an $\loc$ form. 

If $2 deg(\phi) < dim X$ then the same argument as above implies that $\po - \pf$ extends as an $\loc$ form. If $2 deg(\phi) = dim X$ then $\po - \pf = 0$ in $\tb - \sig$ since it is a form of degree higher than the dimension of $\bo$. Hence it extends by zero as an $\loc$ form on $X$.
\end{proof}

\section{Residues}
\setcounter{equation}{0}
The residue is in general a  current supported on the singularities of the bundle map. However, when $\alpha$ is normalized at the singularities the residue is a smooth form.
\begin{lemma}
Let $\alpha : E \rightarrow F$ be an injective map outside $\bigcup \sig$. Suppose that $\alpha$ is normalized at each $\sig$.  Then the residue $\res$  is a  smooth differential form supported on $\sig$ and is given by
$$\res \equiv  \int_{\pi_{i} \mid_{\tb}} T$$
where $\pi_{i} :  \bo \rightarrow \sig$ is projection. Furthermore if either $rank E = rank F$ or $2 deg(\phi) = dim X$ then $\res$ is closed.
\end{lemma}
\begin{proof}
 The normalization condition implies that the transgression form $T$ is pullback form $\bo$. This radial invariance makes $T\mid_{\bo}$ essentially independent of $\epsilon$. In particular
$$\im \int_{\pi_{i} \mid_{\tb}} T = \int_{\pi_{i} \mid_{\tb}} T$$
for any sufficiently small $\epsilon$. Since $T$ is a smooth form outside the singularities, integrating over the fibres of the projection map $\pi_{i}$ yields a smooth closed form on $\sig$.

As mentioned in Remark 5.4 if rank(E) = rank(F) then $\phi(\df) - \phi(\de)$ extends as a closed smooth form on $X$, hence $\res$ is closed. Also, as observed in the proof of Theorem 6.3, the normalization of $\alpha$ implies the $\po -\pf$ is a pullback from $\bo$. If $2 deg(\phi) = dim X$ then $\po - \pf = 0$ in $\tb - \sig$, since it is a form of degree higher than the dimension of $\bo$. So $\po - \pf$ extends by zero to be a closed current and $\res$ is closed.

\end{proof}

\section{Homotopy, Normalized Maps, and Normalized Bundles}
\setcounter{equation}{0}
We showed in \S 6 that $\alpha$, $E$ and $F$ had to be normalized at the singularities $\sig$ for $T$ and $dT$ to extend as $\loc$ forms over the manifold $X$. We now show that normalization is not a strong condition. More precisely we prove that any bundle $E$ with connection is smoothly homotopic to a normalized bundle at $\sig$ and that any bundle map $\alpha$ with singularities $\sig$ is smoothly homotopic to a normalized map on $\cuts$. We also show that any two normalized maps with singularities $\sig$ are homotopic through normalized maps.

\begin{lemma}
Any pair $(E,\de)$ is smoothly homotopic to a normalized bundle at a submanifold $\Sigma$ of $X$.
\end{lemma}
\begin{proof}
 Let $N_{\epsilon/2}$ be a tubular neighborhood of $\Sigma$. Define a map $\pi : X \rightarrow X$ as follows.
$$\pi(v) = \lambda(\norm) v$$
where $\lambda : [0,\infty] \rightarrow [0,1]$ is a smooth function with the properties
$$\lambda(s) = 0 \quad for \quad s<\epsilon/2$$
$$\lambda(s) = 1 \quad for \quad s \geq \epsilon$$
and
$$\lambda '(s) \geq 0$$
Then $\pi \mid_{N_{\epsilon/2}} = p : N_{\epsilon/2} \rightarrow \Sigma$, where $p$ is projection onto $\Sigma$. Hence $\pi^{\ast}(E,\de)$ is equivalent to $E$ and is normalized at $\Sigma$. 

Now let $\pi_{t}(v) = (1-t)v + t\pi(v)$ for $0 \leq t \leq 1$. Define 
$$(E_{t},D_{t}) \equiv \pi_{t}^{\ast} (E,\de)$$
then $E_{t} \cong E$ for $0 \leq t \leq 1$ and $\pi_{1}^{\ast} (E,\de) = \pi^{\ast} (E,\de)$. Hence $\pi_{t}^{\ast}$ is the required homotopy.

\end{proof}

\begin{lemma}
Any bundle map $\alpha : E \rightarrow F$ with singularities $\sig$ is homotopic to a normalized map on $\cuts$
\end{lemma}
\begin{proof}
 Let $\tb$ be an $\epsilon$-tubular neighborhood of $\sig$. We first normalize $E$ and $F$ at each $\sig$. Define a map  $\rho : \cuts \rightarrow \cuts$ as follows. For $v\in \tb -\sig$, set
$$\rho(v) = l(\norm)v$$
where $l : (0,\infty] \rightarrow [1,\infty]$ is a smooth function with the properties 
$$l(t) = \frac{1}{t} \quad for \quad t \leq \epsilon/2$$
$$l(t) = 1 \quad for \quad t \geq \epsilon$$
and
$$l'(t) \leq 0$$
Then $\rho \mid_{N_{i,\epsilon/2}} : N_{i,\epsilon/2} \rightarrow \partial N_{i,\epsilon/2}$ is radial projection in the sense of Definition 6.2, and extends smoothly to be the identity outside $\tb$. Hence $\rho^{\ast} \alpha$ is normalized at $\sig$.
Now let $\rho_{t}(v) = (1-t)v + t\rho(v)$ for $0 \leq t \leq 1$. Then $\rho_{t}^{\ast} (\alpha)$ is the required homotopy between $\alpha$ and $\rho^{\ast}\alpha$ defined on $\cuts$.

\end{proof}

By the same argument as above we have the following lemmas.
\begin{lemma}
Let $\alpha_{1} : E \rightarrow F$ and $\alpha_{2} : E \rightarrow F$ be normalized at $\sig$, with the same tubular neighborhood structure such that
$$\alpha_{1} = \alpha_{2} \on X - \bigcup\tb$$ 
Then there is a homotopy between $\alpha_{1}$ and $\alpha_{2}$ on $\cuts$ through normalized bundle maps.
\end{lemma}

\begin{lemma}
Let $(E_{1},D_{E_{1}})$ and $(E_{2},D_{E_{2}})$ be  normalized at $\sig$, with the same tubular neighborhood structure such that
$$(E_{1},D_{E_{1}}) = (E_{2},D_{E_{2}}) \on X - \bigcup\tb$$
Then there is a homotopy between $(E_{1},D_{E_{1}})$ and $(E_{2},D_{E_{2}})$ on $X$ through  normalized bundles.
\end{lemma}

\section{Invariance of Residues under Homotopy}
\setcounter{equation}{0}

We discuss what happens to the residue when we homotope $\alpha$ through normalized maps at the singularities $\sig$. First we recall a double transgression formula found in [HL1].

\begin{lemma}
Let $D_{s,t}$ be a 2-parameter family of connections, $0 \leq s \leq 1$, and $a \leq t \leq b$, with $D_{s,a} = D_{a}$ and $D_{s,b} = D_{b}$ for all $ 0 \leq s \leq 1$. Then the two transgressions, $T_{1}$ and $T_{0}$, determined by $D_{1,t}$ and $D_{0,t}$ satisfy 
$$T_{1} - T_{0} = dR$$
with
$$R = \int_{a}^{b} \int_{0}^{1} \phi (\frac{\partial}{\partial s} w_{s,t} ; \frac{\partial}{\partial s} w_{s,t} ; \Omega_{s,t}) ds dt$$
where 
$$\phi (A,B;C) = \frac{\partial^{2}}{\partial s \partial t} \phi (C + sA + tB) \mid_{s = t = 0}$$
\end{lemma}

The lemma above allows us to prove the following invariance of the residue classes for the equirank case.

\begin{theorem}
Let $\alpha_{0} : E \rightarrow F$ and $\alpha_{1} : E \rightarrow F$ be normalized maps at singularities $\sig$.  Assume that rank(E) = rank(F). Then the residues, $\res^{0}$ and $\res^{1}$, define the same cohomology class on $\sig$.
\end{theorem}                                                                                   \begin{proof}
Since $E$ and $F$ are of the same rank and the normalization conditions are satisfied we have that
$$\phi(\df) - \phi(\de) = \su \res^{0} [\sig] + dT_{0} \quad\mbox { for }\quad \alpha_{0}$$
and
$$\phi(\df) - \phi(\de) = \su \res^{1} [\sig] + dT_{1} \quad\mbox { for }\quad \alpha_{1}$$
where
$$\res^{0} \equiv  \im \int_{\pi_{i} \mid_{\bo}} T_{0}$$
and
$$\res^{1} \equiv  \im \int_{\pi_{i} \mid_{\bo}} T_{1}$$

By lemma 8.3 we can write a smooth homotopy $\alpha_{t}$ between $\alpha_{0}$ and $\alpha_{1}$ through normalized maps, hence $T_{t}$ extends as an $\loc$ form on $X$ for all $0 \leq t \leq 1$. The initial and end points for the two parameter family of pushforward connections defined by the homotopy are $\df$ and $\de$. Thus we are in the setting of the double transgression lemma above and we can write
$$T_{0} - T_{1} = dR \on  X - \bigcup \sig$$
Since $T_{0}$ and $T_{1}$ extend as $\loc$ forms on $X$ so does $dR$. Furthermore $R$ extends as an $\loc$ form because it is a pullback under radial projection.

We now apply the same argument as we used in Theorem 5.1. Since $T_{0}$ and $T_{1}$ extend as $\loc$ forms then $\displaystyle{\im} (T_{0} - T_{1})\wedge [X-\bigcup \tb] = (T_{0} - T_{1}) \wedge [X]$. Therefore
\begin{eqnarray*}
(T_{0} - T_{1}) \wedge [X] & = & \im dR \wedge [X-\bigcup \tb]\\
                           & = & \im dR \wedge [X- \bigcup\tb] + \im\su R\wedge\bo\\
                           & = & dR\wedge[X] + \im\su R\wedge\bo\\
\end{eqnarray*}
Using the convergence lemma in Federer [F;4.1.19], as before we get that $\displaystyle{\im}R\wedge\bo$ exists. Now we observe that radial invariance makes $R\mid_{\bo}$ essentially independent of $\epsilon$. In particular
\begin{eqnarray*}
\im \int_{\pi_{i}\mid_{\tb}}R & = & \int_{\pi_{i}\mid_{\tb}}R\\
                             & \equiv & R_{0}\\
\end{eqnarray*}
for any sufficiently small $\epsilon$. Therefore $R_{0}$ is a smooth form on $\sig$. This shows that
\begin{eqnarray*}
\res^{0} - \res^{1} & = & \im \int_{\pi_{i}} dR\\
                    & = & d\{\im\int_{\pi_{i}}R\}\\
                    & = & dR_{0}\\
\end{eqnarray*}
Hence they define the same cohomology class of $\sig$.

\end{proof}

\section{The Universal Transgression Form}
\setcounter{equation}{0}
In \S 3 we wrote down the following transgression formula which is valid on $\ho$, the bundle of injective maps from $E$ to $F$.
$$\uni - \unf = d\tilde{T} \on \ho$$
where $\pi : \ho \rightarrow X$ is the projection map. If $\alpha : E \rightarrow F$ is an injective map outside $\sig$ then the equation above pulls down by $\alpha$ to give 
$$\po - \pf = dT \on X - \bigcup \sig$$
In particular, the transgression form $T$ is just a pullback of the universal transgression form $\tilde{T}$ outside the singularities, i.e.,
$$ T = \alpha^{\ast} \tilde{T} \on X - \bigcup \sig$$
This implies that the residues can be expressed universally in terms of $\tilde{T}$ and the map $\alpha : \cuts \rightarrow \ho$, considered as a section of $\ho$ outside the singularities. More precisely,
\begin{equation}
 \res =  \im \int_{\pi_{i} \mid_{\bo}} \alpha^{\ast} \tilde{T}
\end{equation}
In a sense, this is a generalization of the notion of the index of a vector field, in that the residue measures the twisting of the map $\alpha$.

\section{Obstructions}
\setcounter{equation}{0}
The most interesting applications of the main residue formula (5.1) arise when rank(E) = rank(F). For any invariant polynomial $\phi$, the characteristic form $\phi(\df) - \phi(\de)$ extends as a closed, smooth differential form on $X$. Suppose the singularities $\sig$ of the bundle map $\alpha : E \mapsto F$ are \textbf{orientable}, closed submanifolds of $X$. Furthermore, assume that $\alpha$ is extendable. The condition that the $\sig$'s are orientable arises in many natural settings, for instance, orientation preserving finite group actions on a compact manifold. The residues $\res$ are then closed currents on $X$ (see Remark 5.4). More precisely, by Theorem 5.1, $\res \in H^{2deg\phi - e_{i}} (X ; \rb)$, where $e_{i} = codim \res$. We then have the following immediate corollary.

\begin{corollary}
Let $\alpha : E \mapsto F$ be an extendable bundle map with singularities $\sig$ that are orientable, closed submanifolds of $X$. Assume that $rank E  = rank F$. Let $dim X = k$ and $codim \sig = e_{i}$. Suppose that the cohomology groups $H^{2deg\phi - e_{i}} (X ; \rb) = 0$ for each $i$. Then $\phi(\df)$ and $\phi(\de)$ are cohomologous.
\end{corollary}
\begin{proof}
 By (5.1) we have that
$$ \phi(\df) - \phi(\de) = \su \res [\sig] + dT$$
Since $H^{2deg\phi - e_{i}} (X ; \mathbf{R}) = 0$, each $\res$ is exact, i.e.,
$$\res = dS_{i}$$
for some $S_{i} \in \Omega^{2deg\phi - e_{i} - 1} (X ; \rb) $
Hence
$$ \phi(\df) - \phi(\de) = d ( \su S_{i} [\res] + T )$$
which proves the assertion.

\end{proof}
This corollary can be viewed as an obstruction theorem to the existence of bundle maps with orientable singularities of a certain codimension.

\begin{corollary}
Let $E$ and $F$ be vector bundles over $X$ where $rank E = rank F$. Suppose that $\phi(\de)$ and $\phi(\df)$ are not cohomologous for a given invariant polynomial $\phi$. Also suppose that $H^{2deg\phi - k} (X ; \rb) = 0$ for some $k < 2 deg \phi$. Then there cannot exist a bundle map $\alpha : E \mapsto F$ with an orientable, closed singularity $\Sigma$ of codimension $k$.
\end{corollary}
\begin{proof}
 Assume that such an $\alpha$ exists. We can always homotope $\alpha$ such that it becomes extendable with the same singularity $\Sigma$. By Corollary 11.1 $\phi(\de)$ and $\phi(\df)$ are cohomologous, which contradicts the assumptions.

\end{proof}
\nopagebreak
We also the have the following.

\begin{corollary}
Let $\alpha : E \mapsto F$ be an extendable bundle map with singularities $\sig$, where $rank E = rank F$ and $codim \sig = e_{i}$. Let $\phi$ be an invariant polynomial such that $2 deg \phi < max \hspace{2mm}codim \sig$. Then $\phi(\de)$ and $\phi(\df)$ are cohomologous.
\end{corollary}
\begin{proof}
 Again by Theorem 5.1, if  $2 deg \phi < max \hspace{2mm} codim \sig$ then $\res = 0$  for each $i$. Hence
$$\phi(\df) - \phi(\de) = dT$$

\nopagebreak
\end{proof}

\section{Singularities of Maps}

Let $X$ and $Y$ be smooth Riemannian manifolds of equal dimension and consider a smooth mapping
$$f : X \rightarrow Y$$
Consider the differential map 
$$df : TX \rightarrow f^{\ast}TY$$
Here we endow $TX$ and $f^{\ast}TY$ with the standard riemannian connections, normalized at the singularities of $df$. We are now in the standard setting where $E = TX$ and $F = f^{\ast}TY$. Let $p_{k}(Y)$ and $p_{k}(X)$ be the k-th Pontryjagin forms in the normalized riemannian curvatures of $X$ and $Y$.

A straightforward application of the main residue theorem yields the following result.

\begin{theorem}
Suppose that $f : X \rightarrow Y$ is a smooth map between compact oriented riemannian n-manifolds. Suppose that the differential map $df$ has submanifold singularities $\sig$ and  is extendable. Then for any $k \leq n/4$
$$f^{\ast} p_{k}(Y) - p_{k}(X) = \su Res_{p_{k},i} [\sig] + dT$$
and $\res=0$ if $codim\sig >  4k$.
\end{theorem}
\begin{proof}
 We just observe that $p_{k}(f^{\ast}TY) = f^{\ast}p_{k}(TY)$ and apply Theorem 5.1.

\end{proof}

An interesting special case occurs when n = 4k. This yields a 4k-dimensional analogue of the classical Riemann-Hurwitz theorem [M1],[M2],[R].

\begin{corollary}
Suppose that $f : X \rightarrow Y$ is a smooth map between compact oriented 4k-manifolds. Suppose that the differential map $df$ has submanifold singularities
$\sig$ and is extendable. Then
$$M_{f} p_{Y} - p_{X} = \su \int_{\sig} Res_{p_{k},i}  $$
where $M_{f}$ is the degree of the map $f$ and $p_{Y}$ and $p_{X}$ are the top Pontryjagin numbers of the manifolds $X$ and $Y$ respectively.
\end{corollary}

More generally we consider $\wp$, a homogeneous polynomial of weight k in k indeterminates. The associated Pontryjagin number to $\wp$ of a compact, oriented manifold $X$ of dimension 4k is defined to be
$$\wp(X) = \int_{X} \wp(p_{1}(X), \ldots , p_{k}(X))$$
Then we have the following corollary.
\begin{corollary}
Suppose that $f : X \rightarrow Y$ is a smooth map between compact oriented 4k-manifolds. Suppose that the differential map $df$ has submanifold singularities
$\sig$ and   is extendable. Then
$$M_{f} \wp(Y) - \wp(X) = \su \int_{\sig} Res_{\wp,i} $$
where $M_{f}$ is the degree of the map $f$ and $\wp(X)$ and $\wp(Y)$ are the               Pontryjagin numbers associated to $\wp$ of the manifolds $X$ and $Y$ respectively .
\end{corollary}

For a topological approach to these formulae see [N], [GGV].

In particular, by the Hirzebruch signature formula [MS] which states that the signature of a compact, oriented manifold of dimension 4k is expressible as a polynomial in the Pontyjagin classes through the L-class, we have the following.

\begin{corollary}
Suppose that $f : X \rightarrow Y$ is a smooth map between compact oriented 4k-manifolds. Suppose that the differential map $df$ has submanifold singularities
$\sig$ and   is extendable. Then
$$M_{f} sig(Y) - sig(X) = \su \int_{\sig} Res_{sig,i} $$
where $M_{f}$ is the degree of the map $f$ and sig(X) and sig(Y) are the               signatures of the manifolds $X$ and $Y$ respectively.
\end{corollary}

\begin{remark}
The results in this section apply naturally to branched coverings
$$ \pi : X \mapsto Y$$
branched along a submanifold $\Sigma$ of codimension 2 in $Y$ and also to finite group actions on a compact manifold. In a subsequent paper we study these important cases where the singularities are of nongeneric codimension in more detail, providing explicit calculations of the residues that appear in the formulae above.
\end{remark}

\begin{remark}
It is clear that similar formulae hold for the case of maps between complex manifolds.
\end{remark}

\section{CR-Singularities}
The methods introduced above yield interesting results in CR-geometry. Consider an immersion 
$$f : X \hookrightarrow Z$$
of a real manifold $X$ into a complex manifold $Z$ where,
$$ dim_{\rb}(X) = n = dim_{\cb}(Z)$$
Then the differential map 
$$df : TX \rightarrow f^{\ast}TZ$$
extends to a complex bundle map
$$df_{\cb} : TX \otimes_{\rb} \cb \rightarrow f^{\ast}TZ$$

Assume that the bundle map $df_{\cb}$ has submanifold singularities $\sig$. This corresponds to the loci of points where $f_{\ast} T_{x}X$ contains a complex subspace having `excess' dimensions, i.e., more complex tangency than expected. Specifically we have:

\begin{lemma}
$$\bigcup \sig = \{x \in X : dim_{\cb}(T_{x} \cap JT_{x}X) > 0\}$$
where $J$ denotes the complex structure of $Z$.
\end{lemma}
\begin{proof}
 We have that $\bigcup \sig = \{x \in X : rank(df_{\cb}) < n\} $. Note that at $x \in X$,
\begin{eqnarray*}
ker(df_{\cb}) & = & \{V + iW : V,W \in T_{x}X\quad\mbox{and}\quad V + JW = 0\}\\
              & = & \{V + iJV : V,JV \in T_{x}X\}\\
              & = & [(T_{x}X \cap JT_{x}X) \otimes \cb]^{0,1}\\
\end{eqnarray*}
Since rank($df_{\cb}$) = n - $dim_{\cb}ker(df_{\cb})$ we have rank($df_{\cb}$) $<$ n if and only if $dim_{\cb}(T_{x}X \cap JT_{x}X) > 0$.

\end{proof}

We now suppose that $X$ carries a Riemannian metric and $Z$ carries a hermitian metric and a complex connection, and we normalize these connections at each $\sig$. Let $p_{i}(X)$ and $c_{i}(Z)$ be the ith Pontryjagin and Chern forms of $X$ and $Z$ respectively.

Applying Theorem 5.1 yields the following result.

\begin{theorem}
Let $f : X \hookrightarrow Z$ be an immersion of a real m-manifold into a complex m-manifold. Assume that $df_{\cb}$ has submanifold singularities $\sig$ and that it is extendable. Then
$$f^{\ast} c_{2i}(Z) - p_{i}(X) = \su \res [\sig] + dT$$
\end{theorem}
\begin{proof}
 Apply Theorem 5.1 to the bundle map $df_{\cb}$ and observe that
$$(-1)^{i}c_{2i}(TX \otimes \cb) = p_{i}(X)$$

\end{proof}

Consider the case when $Z = \cb^{n}$. This gives the following.

\begin{corollary}
Let $f : X^{n} \hookrightarrow \cb^{n}$ be an immersion with the property that $df_{\cb}$ has submanifold singularities $\sig$ and that it is extendable. Then
$$p_{i}(X) = \su \res [\sig] + dT$$
\end{corollary}

\section{Finite Singularities and a Generalized Hopf Index Formula}
The finite singularity sets of r-vector fields and r-plane fields were studied extensively by E.Thomas [T1],[T2], Atiyah [A], and Atiyah-Dupont [AD].
In certain cases they managed to relate these singularities to algebraic invariants of the manifold. We continue the study of finite singularities in the general setting of bundle maps over a compact manifold.

In \S 9 we discussed the universal transgression form $\tilde{T}$. We now use the universal construction to prove the following theorem.

\begin{theorem}
Let $\alpha : E \rightarrow F$ be a bundle map over a compact manifold $X$ with isolated finite singularities $\{x_{i}\}$. Assume that $\alpha$ is normalized at each $x_{i}$. Let rank(E) = m and rank(F) = n. Then for any invariant polynomial of top degree,
$$\int_{X} \po-\pf=\su \int_{S^{i}} \alpha^{\ast}\tilde{T}$$
where $\tilde{T}\in\Omega^{odd}(\stief)$. Here $\stief$ is the Stiefel manifold of m-frames in $\rb^{n}$ and the $S^{i}$ are small spheres around the points $\{x_{i}\}$.
\end{theorem}
\begin{proof}
 We recall that $T=\alpha^{\ast}\tilde{T}$ outside $x_{i}$. Since $E$ and $F$ are normalized at $x_{i}$ then near $x_{i}$
$$\ho = Hom^{\ast}(\rb^{m},\rb^{n}) = \stief$$
Now we apply Theorem 5.1 and integrate over the manifold. Since $\phi$ is of top degree by Lemma 7.1 we do not need to take a limit for the residue.

\end{proof}

This theorem is a bundle map analogue of the classical Hopf index formula.
It has the following  interesting corollary.
\begin{corollary}
Consider a map $\alpha : E \rightarrow F$ over a compact manifold $X$ of dimension 4n with isolated finite singularities, where rank(E) = 2 and rank(F) = n. Suppose that $\alpha$ is normalized at the singularities. Then
$$\int_{X} p_{n}(F)-p_{1}(E)p_{n-1}(I^{\perp}) = 0$$
In particular if $p_{1}(E) = 0$ then $p_{n}(F)=0$. Here $p_{i}$ denotes the i-th Pontryjagin class.
\end{corollary}
\begin{proof}
This is just a consequence of dimension. By the theorem above $\tilde{T}$ is a differential form of degree 4n - 1 on $V_{4n,2}$. But $dim V_{4n,2} = 4n - 3$ and hence $\tilde{T}=0$. 

\end{proof}

\begin{remark}
The result above is not a consequence of obstruction theory because in general, $\pi_{4n-1}(V_{4n,2}) \neq 0$.
\end{remark}

We would like to know when $\tilde{T}$ is closed near $x_{i}$. Then the residue
can be interpreted using cohomology of $\stief$. For this we use the following construction. 

Let $\gras$ be the Grassmann manifold of m-planes in $\rb^{n}$ and let $\rho : \stief \rightarrow \gras$ be the standard fiber map. Also let $\tau$ be the tautological bundle over $\gras$ and $\tau^{\perp}$ be its dual. Give $\tau^{\perp}$ the connection induced by projection. This construction yields the following lemma.

\begin{lemma}
$\tilde{T}$ is closed near $x_{i}$ iff $\phi(\tau^{\perp}) = 0$ on $\gras$.
\end{lemma}
\begin{proof}
 We can write
$$\phi(\rb^{n}) - \phi(\rb^{m}\oplus\tau^{\perp}) = d\hat{T} \on \gras$$
By construction $\tilde{T}=\rho^{\ast}\hat{T}$. The equation above reduces to
$$\phi(\tau^{\perp})=d\hat{T}$$
Therefore if $\phi(\tau^{\perp}) = 0$ on $\gras$, $\hat{T}$ is closed and hence $\tilde{T}$ is also closed.

\end{proof}

The lemma above yields the following corollary.
\begin{corollary}
Let $\alpha : E \rightarrow F$ be a bundle map over a compact manifold $X$ of dimension 4i with isolated finite singularities $\{x_{i}\}$. Assume that $\alpha$ is normalized at each $x_{i}$. Let rank(E) = m and rank(F) = n. Suppose that $2i>n-m$. Then
$$\int_{X} p_{i}(F) - p_{i}(E\oplus I^{\perp}) = \su \int_{S^{i}} \alpha^{\ast}\tilde{T}$$
where $\tilde{T} \in H^{4i-1}(\stief)$.
\end{corollary}
\begin{proof}
 We observe that $rank(\tau^{\perp}) = n-m$. If $2i>n-m$ then
$$0=p_{i}(\tau^{\perp}) = d\hat{T}$$
Hence $\hat{T}$ is closed and $\tilde{T} \in H^{4i-1}(\stief)$.

\end{proof}

\begin{corollary}
Let $\alpha : E \rightarrow F$ be a bundle map over a compact manifold $X$ of dimension 2i with isolated finite singularities $\{x_{i}\}$. Assume that $\alpha$ is normalized at each $x_{i}$.Let rank(E) = m and rank(F) = n. Let $\phi$ be any multiplicative invariant polynomial and $\phi_{i}$ be its i-th degree term. Then for $m(n-m)<2i\leq m(n-m) + \frac{1}{2}m(m-1)$ we have that
$$\int_{X} \phi_{i}(F) - \phi_{i}(E\oplus I^{\perp}) = \su \int_{S^{i}} \alpha^{\ast}\tilde{T}$$
where $\tilde{T} \in H^{2i-1}(\stief)$
\end{corollary}
\begin{proof}
 This again is a consequence of dimension. We have that 
$$dim \stief = m(n-m) + \frac{1}{2}m(m-1)$$
and
$$dim \gras = m(n-m)$$
If $2i>dim\gras$ then
$$0=\phi(\tau^{\perp})=d\hat{T}$$
because it is a differential form of degree higher than the dimension of $\gras$. Therefore $\hat{T}$ is closed and hence $\tilde{T} \in H^{2i-1}(\stief)$.

\end{proof}

\section{Clifford and Spin Bundles}
Let $\pi : F \mapsto X$ be a 2n-dimensional vector bundle with spin structure (for a complete discussion of the following constructions see [LM]). Assume that $F$ is provided with a connection $D$. Let $ \mathcal{S}$ denote the complex spinor bundle associated to $F$ and let $D_{\mathcal{S}}$ be the connection on $\mathcal{S}$ induced from the one on $F$. There is a canonical decomposition 
$$ \mathcal{S} = \mathcal{S^{+}} \oplus \mathcal{S^{-}}$$
by the complex volume form. Suppose that we are given an odd form $\alpha$, i.e., $\alpha \in \Gamma (\bigwedge^{odd} (F) )$. Then $\alpha$ is a bundle map from $\mathcal{S^{+}}$ to  $\mathcal{S^{-}}$. Assume that
 $$\alpha : \mathcal{S^{+}} \mapsto \mathcal{S^{-}}$$
has closed submanifold singularities $\sig$. Then again we are in the setting of theorem 5.1, with rank $\mathcal{S^{+}}$ = rank $\mathcal{S^{-}}$ . Furthermore , if $E$ is any complex vector bundle with complex vector bundle with complex connection $D_{E}$, we let 
$$D_{\mathcal{S} \otimes E} = D_{\mathcal{S}} \otimes 1 + 1 \otimes D_{E}$$
denote the tensor product connection on $\mathcal{S} \otimes E$ and normalize this connection at $\sig$. This connection observes the splitting of the spinor bundle. The odd form $\alpha$ is again a bundle map from $\mathcal{S^{+}} \otimes E$ to $\mathcal{S^{-}} \otimes E$ and we have the following.

\begin{corollary}
Let $(F , D)$ and $(\mathcal{S} \otimes E , D_{\mathcal{S} \otimes E})$ be as above. Suppose that $\alpha \in \Gamma (\bigwedge^{odd} (F) )$ has closed submanifold singularities  $\sig$ as a bundle map $\alpha : \mathcal{S^{+}} \mapsto \mathcal{S^{-}}$ and is extendable. Then
$$ ch (D_{\mathcal{S^{+}} \otimes E}) - ch(D_{\mathcal{S^{-}} \otimes E}) = \su Res_{ch,i} [\sig] + dT$$
\end{corollary}

\begin{remark}
The same formula holds when $F$ is $spin^{c}$.
\end{remark}

Even if $F$ is not spin or $spin^{c}$, we can derive an interesting residue formula by considering the complex Clifford bundle $C\ell \equiv C\ell \otimes \mathbf{C}$ associated to $F$. Again we have a decomposition
$$ C\ell = C\ell^{+} \oplus C\ell^{-}$$
and an odd form $\alpha \in \Gamma (\bigwedge^{odd} (F) )$ acts as a bundle map
$$\alpha : C\ell^{+} \otimes E \mapsto C\ell^{-} \otimes E$$
where $E$ is any complex vector bundle with complex connection $\de$. Again we let
$$D_{C\ell \otimes E} = D_{C\ell} \otimes 1 + 1 \otimes D_{E}$$
denote the tensor product connection on $C\ell \otimes E$ and we normalize this connection at $\sig$, the submanifold singularities of $\alpha$. As before this connection preserves the splitting of $C\ell$. We immediately have

\begin{corollary}
Let $(F,D)$ and $(C\ell \otimes E, D_{C\ell \otimes E})$ be as above . Suppose that $\alpha \in \Gamma (\bigwedge^{odd} (F) )$ has closed submanifold singularities as a bundle map $\alpha : C\ell^{+} \otimes E \mapsto C\ell^{-} \otimes E$ and is extendable. Then
$$ ch (D_{C\ell^{+} \otimes E}) - ch(D_{C\ell^{-} \otimes E}) = \su Res_{ch,i} [\sig] + dT$$
\end{corollary}

\begin{remark}
An explicit calculation of the residue in these two cases would be useful because it would yield analogues of Grothendieck - Riemann - Roch [AH].
\end{remark}

\end{document}